\documentstyle[a4,twoside]{article}
\makeatletter
\renewcommand{\@oddhead}{Operatoren ohne Dichotomie  \hfill \thepage}
\renewcommand{\@evenhead}{\thepage \hfill Sergej A. Choro\v savin }
\renewcommand{\@oddfoot}{}
\renewcommand{\@evenfoot}{}
\makeatother

{\vspace*{2ex}\par {#1}\quad\it }{\medskip\par }

\newenvironment{Thm}[2]%
{\par\addvspace{\bigskipamount}{\bf #1#2}\it }{\par\addvspace{\bigskipamount} }

\newenvironment{Definition}[1]{\begin{Thm}{Definition}{#1}}{\end{Thm} }

\newenvironment{Satz}[1]{\begin{Thm}{Satz}{#1}}{\end{Thm} }
\newenvironment{Hilfssatz}[1]{\begin{Thm}{Hilfssatz}{#1}}{\end{Thm} }

\newenvironment{Observation}[1]{\begin{Thm}{Observation}{#1}}{\end{Thm} }
\newenvironment{Bemerkung}[1]{\begin{Thm}{Bemerkung}{#1}}{\end{Thm} }

\newenvironment{Beweis}{\par\addvspace{\bigskipamount} {\sc Beweis}}%
{\par\hspace*{\fill}$\Box$\par\addvspace{\bigskipamount} }

\newcommand{\dfrac}[2]{{\displaystyle \frac{#1}{#2}}}

\begin{document}
\begin{center}
{\bf Hamiltonsche Bahnen ohne  Zerspaltungseigenschaft.
 Die L\"osung einer Aufgabe von  M.G. Krein.}
\bigskip

 {{\sc Von  Sergej A. Choro\v savin  } }
\bigskip

\bigskip

 Email: { sergius@pve.vsu.ru }

\end{center}

\begin{abstract}
 Here we construct linear hamilton systems without usual dichotomy
 property. The Ljapunov spectra of these systems and
 the behaviour of trajectories are very complicated. 
 The paper's subject refers to some problems 
 of indefinite inner product methods
 in the stability theory of abstract dynamical equation solutions .
\end{abstract}

\section%
 { Einleitung}

 Man betrachte Bahnen irgendeines dynamischen linearen hamiltonschen Systems.
 Gibt es welche Bahnen, die vorgegebenes 
 Wachstumsverhalten f\"ur $t\to \pm\infty$ haben ?
 Diese Frage sowie die Frage, ob die ensprechenden 
 Zerlegungen von beliebigen Bahnen vorhanden sind,
 diese beiden Fragen sind traditionell und
 immer noch zeitgem\"a{\ss} f\"ur Experten in
 Theorie von dynamischen Systeme.
 Noch mehr, 
 wie man diese Fragen beantwortet,
 davon h\"angt die Gestalt und der Inhalt von vielen Seiten
 von der analytische Stabilit\"atstheorie ab.

 Solange man sich mit einem linearen endlich-dimensionalen
 Hamiltonschen Systeme
 besch\"aftigt, 
 hat man die wohlbekannte Antwort:

\medskip

{\it
 Jede Bahn $x(t)$ l\"a{\ss}t sich auf folgende Weise darstellen:
$$
 x(t)=x_1(t)+\cdots+x_N(t)    \leqno(*)
$$
 wobei  $x_1(t),\ldots,x_N(t)$  die Bahnen des urspr\"unglichen Systems sind,
 die f\"ur \hbox{$t\to\infty$} das 
\"ubliche
 exponentiellweise Wachstumsverhalten haben.
}

 Insbesondere gilt

\medskip

 {\bf Zerspaltungssatz 1.}{\it
 Beliebige Bahn $x(t)$ wird folgenderma{\ss}en dargestellt:
$$
 x(t)=x_-(t)+x_+(t)\,;\,|x_\pm(t)|\leq P_\pm(t) 
 \hbox{ \rm f\"ur } t\to\pm\infty
\leqno(**)\,
$$
 {\rm(} hier sind $P_\pm$ gewisse Polynome und $x_\pm$ sind geeignete Bahnen  
 desselben dynamichen Systems  {\rm)}
}.

\medskip

 Soweit ist alles einfach.
 {\bf Wie steht es aber mit den Systemen,
 die unendlich-dimensional sind?}

 Es ist dabei ersichtlich, da{\ss} die obige Zerspaltungseigenschaft
 umformuliert werden mu{\ss}.
 Mindestens zwei Umformulierungen scheinen geeignet zu sein:


\medskip\noindent  
 1)\hfill 
\parbox[t]{0.95\textwidth}{
statt in $(**)$ eingetretenen Polynome 
 habe man
 ``nicht zu schnell
 wachsende''  \\n\"ahmlich die ``subexponentiellwachsende'' Funktionen
 zu gestatten;
 dabei sprechen  wir \"uber die Zerspaltung von Bahnen 
 und stellen wir die zugeh\"orige {\bf Existenzaufgabe} von solchen Bahnen;
}

\medskip\noindent 
 2)\hfill 
\parbox[t]{0.95\textwidth}{
 statt die Menge von subexponentiellwachsenden Bahnen zu suchen,
 habe man 
 diejenige invariante Teilr\"aume
 des Propagators des Systems
 zu suchen, 
 auf denen der
 Spectralradius dieses Propagators
$\leq 1$ sei;
 hierbei spricht man \"uber die {\bf Aufgabe von {\sc M.G. Krein}}.
}

\medskip\noindent 
 Hier ist zu beachten, da{\ss} nicht nur die kontinuierliche ``Zeit''
 (d.h. $t\in {\bf R}$),
 sondern auch die diskrete ``Zeit''
 (d.h. $t\in {\bf Z}$)
 oder irgendeine andere abstrakte ``Zeit''
 von Interesse ist:
 es geht jetzt um ``die symplektischen Darstellungen einer Halbgruppe''.

 Was die Phasenr\"aume, da{\ss} hei{\ss}t
 die R\"aume,
 worauf solche Darstellungen wirken,
 betrifft, 
 mu{\ss} man nicht nur reelle, sondern auch komplexe R\"aume betrachten.
 Es handelt sich um die ``Darstellungen in R\"aumen mit 
 indefiniten inneren Produkt'' sowie um  ``J-unit\"are und J-isometrische
 Darstellungen'' u.\"a.
 Im allgemeinen die Aussage ``$T$ ist $J$-unit\"ar ''
 bedeutet:
$$T^*JT  =  J  =  TJT^*;\quad J^*=J;\quad J^2=I.$$

\underline{Zu erinnern ist an:}

\begin{Hilfssatz}{ H1-1}
$V$ sei ein auf einem Hilbertraum $H$
 definierter linearer beschr\"ankter Operator, der ein
 beschr\"anktes Inverses besitzt.

 Dann ist der Operator $V\oplus V^{*-1}$ (d.h. Hilbertsche 
 direkte Summe des Operators  $V$ mit $V^{*-1}$) 
$J$-unit\"ar in bezug auf  $J$, der durch die folgende 
 Formel erkl\"art ist:
$$
 J\,:\,x\oplus y\,\mapsto\,y\oplus x\quad (x\,,y\, \in H\,) \qquad
$$ 

 Setzt man statt obigem $J$ einen so wirkenden Operator:
$$
 {\cal J}\,:\,x\oplus y\,\mapsto\,y\oplus -x\quad (x\,,y\, \in H\,) \qquad ,
$$
 so erhaltet man, da{\ss}  
${\cal J}^*=-{\cal J};{\cal J}^2=-I$ , d.h. ${\cal J}$ 
 ein Operator von einer symplektischen Struktur ist;
 in diesem Falle ist der Operator
$V\oplus V^{*-1}$ 
 ein symplektischer Automorphismus. \qquad $\Box$

\end{Hilfssatz}

 Wir werden die Konstruktion 
$V\oplus V^{*-1}$ 
 systematisch benutzen und schreiben dann:
${\hat V} := V\oplus V^{*-1}$. 

\medskip

 Um die beiden obigen Aufgaben zu pr\"azisieren, 
 f\"uhren wir eine geeignete Definition ein.

\begin{Definition}{ D1-1}
$$
 S_0(T)\,:=\,\{x\in H|\quad\|T^Nx\|\to 0 \hbox{ f\"ur }N\to +\infty \},
$$ 
$$
 S(T)\,:=\,\{x\in H|\, \exists C\geq 0 \, \forall N\geq 0 
 \quad\|T^Nx\|\leq C \},
$$ 
$$
 S_+(T)\,:=\,\{x\in H|\forall a>1 \, \exists C\geq 0 \, \forall N\geq 0 
 \quad\|T^Nx\|\leq Ca^N\},
$$ 
$$
    r(T):= \hbox{Spektralradius von }T \,.
$$
\end{Definition}

\begin{Bemerkung}{ B1-1}
 Ist $T$ ein linearer beschr\"ankter Operator und ist $L$ ein
$T$-invarianter Teilraum, so da{\ss}  $r(T|L)\leq c$, 
 so ist 
$$
L\subset S_+(c^{-1}T)
$$
\end{Bemerkung}

\begin{Bemerkung}{ B1-2}
$$
 S_{\rm x}(T_1\oplus T_2)=S_{\rm x}(T_1)\oplus S_{\rm x}(T_2);
$$
 hier steht jeweils $S_{\rm x}$ f\"ur $S_0$ oder $S$ oder $S_+$ .
\end{Bemerkung}

\begin{Bemerkung}{ B1-3}
$$
 S_0(T) \perp S(T^{*-1})
 \,;\qquad
 S(T) \perp S_0(T^{*-1})
$$
 Dies folgt schon daraus, da\ss\
$$
 |(x,y)|=
 |(T^Nx,T^{*-N}y)|\leq
 \|T^Nx\|\,\|T^{*-N}y\|   \,.
$$
\end{Bemerkung}

 So besteht die erste Aufgabe darin, die Struktur von der Menge des Typs
$S_{\rm x}(T)$
 zu erkl\"aren.
 Mindestens fragt man:
\begin{Thm}{Frage}{ 1.}

$T$ 
 sei beliebig ausgew\"ahlt, dieser sei aber symplektisch oder $J$-unit\"ar .

 Ist 
$S_+(T) \not= \{0\}$ ?
\end{Thm}

    Die zweite Aufgabe sieht mehr klassisch aus: 
 man suche diejenigen
$T$-invarianten
 Teilr\"aume, 
 worauf die vorgegebene Begrenzung des Spektralradius erf\"ullt ist. 
    Mindestens fragt man:
\begin{Thm}{Frage}{ 2.}

$T$ 
 sei beliebig ausgew\"ahlt, dieser sei aber symplektisch oder $J$-unit\"ar .

 Gibt es einen nicht trivialen $T$-invarianten Teilraum $L$, so da{\ss} 
$r(T|L) \leq 1$ \nolinebreak ?
\end{Thm}

 Diese beiden Fragen sind nicht voneinander unabh\"angig.
 N\"amlich, wir haben (wegen Bemerkung {\bf B1-1}):
 die negative Antwort auf die erste Frage zieht 
 die negative Antwort auf die zweite Frage nach sich
 und die positive Antwort auf die zweite Frage 
 hat zur Folge, da{\ss} die Antwort auf die erste Frage 
 ebenfalls positiv ist.

   Obwohl die beiden obigen Aufgaben vor langem bekannt sind
 (sieh z.B. \cite{Kre64}, [Kre65], [DalKre]), 
 ist die `allgemeine' L\"osung noch immer nicht bekommen.
 Selbst die Antwort auf die Fragen 1 und 2 zu finden,
 das ist erst vor kurzem gelungen.
 Diese Antwort basiert  auf [Ch97], [Ch98]; 
 sie wird hier vorgestellt.


   Man fragt sich zuletzt, ob diese hier 
 er\"orterten
 Aufgaben 
 einigerma{\ss}en 
 beachtenswert 
 sind.

   Herkunft ist mir nicht ganz klar. Ich glaube,
 die Geschichte dieser Aufgaben
 begann am mindesten seit der Zeit, 
 als man den Begriff von Ljapunovschem Index
 \footnote{ 
 man nennt diese Gr\"o{\ss}e auch 
 oberen Wachstumsindex oder Ljapunovsexponent
            }
 eingef\"uhrt hatte
 und lineare Hamiltonsche Gleichungen als ein einzelnen
 merklichen
 Gegenstand betrachtet.

 Es ist vielleicht nicht \"uberfl\"ussig zu erinnern,
 da{\ss} die Ljapunovschen Indizes 
 von einer vorgegebenen Bahn 
$x(t), \quad t \in {\bf R}$
 ( eines linearen dynamischen Systems) 
 reelle Zahlen 
$\lambda_{\pm}$ 
 sind, die durch

$$
  \lambda_{\pm}:=\limsup_{t\to \pm\infty}
                 \frac{\ln \|x(t)\|}{|t|}
  =\inf\{
  \lambda|\: \|x(t)\|\leq C_{\lambda}e^{\lambda |t|}\, (t\to\pm\infty)
  \}
$$ 
 definiert werden.

 So betreffen die oben formulierten Probleme die zugeh\"orige
 Theorie von Ljapunovschem Spektrum, 
 den man manchmal auch Floquet-Spektrum nennt.

 Die erste Spur von vollige Systematisierung und mithin von einer Einf\"uhrung
 von Kanonischen Formen (von endlich-dimensionalen) Hamiltonianen, 
 habe ich in [Arn] als sogennanten Williamsonssatz gefunden
 (originelle Werke sind [Will36] und [Will37]).

 Man beachte: dieser Satz ist in Rahmen von reellen R\"aumen formuliert.

 Es ist kaum m\"oglich, exakt und bestimmt zu erkennen, wer der Erster war,
 der den Phasenraum des linearen Systems formalweise komplexifiziert hatte,
 und 
 statt der symplektischen Struktur (${\cal J}^2=-I$) 
 die Kreinsche $J$-Sturktur ($J^2=I$) 
 eingef\"uhrt hatte.

 Wie dem auch sei, solcher Trick 
 fand sich so gelungen 
 [Kre65], [KL1], [KL2], [IK], [IKL], [DalKre], [DadKul], [Kul], [Ikr89], 
 da{\ss} der besondere Zweig von der Theorie der linearen
 R\"aume und Operatoren---die Theorie von Krein-R\"aume---
 entstand [Bogn], [DR].
\footnote{ Man nennt Krein-R\"aume auch
 "R\"aume mit regul\"arem innerenem Produkte"}


 Der andere Zweig von Entwicklung solchartiger Theorien 
 h\"angt mit solchen Begriffe zusammen, 
 wie {\it lineare kanonische Transfomation},
 {\it kanonische Vertauschunsbeziehungen }, alias
 {\it CCR (Canonical Commutation Relations)},
 {\it Bogoliubov Transformationen} 
 {\it quasifreie Bewegung } 
 [Bog], [BraRob], [Ber], [Rob], [MV], [Emch], [RS2], [RS3], [Oks], 
 [DadKul], [Kul], [Ikr89].
 Die soeben
 erw\"ahnten Begriffe geh\"oren eher zu der 
 Quanten Physik, als zu der reinen Mathematik 
 und, 
 obwohl
 die vorliegende Arbeit vorerst gerade aus diesem Zweig hervorgeht
 [Ch81], [Ch83], [Ch83T], [Ch84], [ChDTh],
 wollen wir auf diese Begriffe nicht eingehen,     
 weil hier nur "rein mathematischer" Teil von unseren Untersuchungen
 dargestellt worden wird.

 Wir sollen noch ein Thema nennen, n\"amlich ``Modelle 
 von chaotischem Verhalten'', oder einfach ``Modelle von Chaos''.  

 Eine verbreitete \"Uberzeugung ist, da{\ss} kein lineares System 
 zu Modellierung des Chaos pa{\ss}t und geeignet ist. 
 Diese \"Uberzeugung ist verbreitet, aber umstritten. 
 Hier meinen wir vorerst 
%
"Does quantum chaos exist? (A quantum Lyapunov exponents approach.)"
 von Wladyslaw Adam Majewski. 
 Wir vertreten einen verwandten Standpunkt und suchen explizite 
 Beispiele von linearen `chaotischen' Systemen zu konstruiren. 

 In der vorliegenden Arbeit 
 konstruiren wir drei linearen diskreten dynamischen Systeme 
 ( jeweilige Operatoren 
$\hat{U}, \hat{V}, \hat{W} $ in der Sektionen 2, 3, )
 und beschreiben kurz ihre kontinuierlichen Analoga 
 (Sektion 4),
 die ein sehr kompliziertes Verhalten haben.

 Es ist nat\"urlich kaum so, 
 da{\ss} man diese Systeme   
 als vollkommen chaotische Systeme 
 qualifizieren wolle. 

 Dennoch sehen ihre Ljapunov Spekren sehr exotisch aus,
 und entschprechen 
 die spektrale Teilr\"aume und die Mengen 
$S_0, S, S_+ $  
 einander etwas seltsam und \"uberraschend.
 Von diesem Standpukte aus, 
 wir w\"urden sagen, das Verhalten dieser Systeme 
 ist {\bf vorchaotisch}.
 Die Situation
 ist
 im gro{\ss}en und ganzen 
 diese:

 Das erste System, das wir konstruieren, ist ein solches, da{\ss} 
$S_+$ 
 ein blo{\ss}es Nullelement enth\"alt. 
 Folglich, alle Ljapunov Indizes sind streng positiv. 
 Noch mehr, w\"ahlt man eine Zahl
$c>1$ 
 ganz beliebig, so kann man ein System konstruiren, 
 so da{\ss} {\bf alle} Ljapunov Indizes 
$>\ln c$.
 Demungeachtet gibt es von Null verschiedene Bahnen, so da{\ss}
 der {\bf untere} Wachstumsindex 
$\leq -\ln c$
 ist.

 F\"ur den Fall des zweiten Systems ist die Menge 
$S_+$ 
 ``reich'': 
 es gibt ``viele'' Bahnen, deren Ljapunov Indizes streng negativ sind
($< -\ln 2 - \ln c$, $c>1$), 
 dann aber hat die Abschlie{\ss}ung dieser Menge 
 gewisse Bahnen, deren Indizes streng positiv sind 
 ($\geq \ln 2 + \ln c$).   

 Was das dritte System betrifft, dieses ist, da{\ss} 
$S_0$ 
 auch von Null verschieden ist 
 ( nat\"urlich, $S_0 \subset S \subset S_+$), 
 sowie ein gesuchter Teilraum existiert, ---
 ein {\bf maximaler} invarianter Teilraum
$L$ derart,
 da{\ss} 
$r(\hat{W}|L) = 1 $
 (folglich, $L \subset S_+$). 
 Doch sind diese, $L$ und $S_0$, zueinander orthogonal und au{\ss}erdem 
$L\cap \overline{S} = \{0\}$.

 Noch mehr,  $L$ =${S_0}^{\perp}$, sogar $L$ =${S}^{\perp}$. 
 Noch mehr , der Spektralradius des auf 
 $L^\perp \equiv \overline{S_0}$ 
 eingeschrenkten Propagators ist gleich 
 $2$; das Spektrum der Einschrenkung selbst liegt in der Menge 
$\{z| 1\leq |z| \leq 2\}$.

 Soweit sind unsere Ergebnisse, die die Frage 1 betreffen. 

 Was die Frage 2 selbst betrifft, 
 die lange Zeit existierte die Ahnunng,
 da{\ss} die Antwort 
%
%
 allemal
%
%
 positiv ist.
 Diese Ahnung unterst\"utzten Ergebnisse von verschiedener Art, 
 sowohl einzelnen Existenzs\"atze --
 von [Kre64] bis [Shk99]--,
 als auch die Theoreme von der Art --
 "Die Menge 
 $S_0(T)$
 liegt in dem Teilraum, der M.G.Krein konstruiert hat [Kre64], [Kre65]"
 (dar\"uber sieh [Ch89-1], [Ch89-2], [Ch96T], [Ch00]).


  Nun 
  gehen wir eigentlich zur Grunddarlegung selbst \"uber.  
  Dabei entspricht die sogenannte `$J$-Terminologie ' der Arbeit
  [Kre65], w\"ahrend die allgemeine mathematische Terminologie 
  auf [RS] folgt.



\newpage
\section
{ Diskretes dynamisches System ohne Dichotomie}

\begin{Satz}{ S2-1}
 Wie auch eine positive Zahl $c>0$ ausgew\"ahlt sein mag,
 es gibt einen  $J$-unit\"aren und symplektischen Operator $\cal U$, da{\ss} 
$$
 S_+(c^{-1}{\cal U})  =  S_+(c^{-1}{\cal U}^{-1})  =  \{0\}\,.
$$

 Insbesondere,

\par\addvspace{\medskipamount}\noindent  
{\rm (i)}\hfill 
\parbox[t]{.93\textwidth}{ 
 $L$ sei ein von Null verschiedener $\cal U$-invarianter Teilraum. 
 Dann ist \mbox{$r({\cal U}|L) > c$; }
}
\par\addvspace{\smallskipamount}\noindent      
{\rm (ii)}\hfill 
\parbox[t]{.93\textwidth}{ 
 $L'$ sei ein von Null verschiedener ${\cal U}^{-1}$-invarianter Teilraum.
 Dann ist \mbox{$r({\cal U}^{-1}|L) > c$; }
}
\par\addvspace{\smallskipamount}\noindent       
{\rm (iii)}\hfill 
\parbox[t]{.93\textwidth}{ 
 es gibt keinen $\cal U$-invarianten Teilraum $L''$,
 so da{\ss}
 \mbox{$|spectrum\,({\cal U}|L'')|\geq c^{-1}$ } .
}
\par\addvspace{\medskipamount}\noindent  


\end{Satz}

\begin{Beweis}{ }

 Angenommen, es existiere ein Operator 
$U$ derart,
 da{\ss} gilt:

\par\addvspace{\medskipamount}\noindent       
{\rm 1)} 
$U$ 
 ist linear, bijektiv, beschr\"ankt;
\footnote{ daraus folgt, da{\ss}  
$U^{*}$, 
$U^{-1}$,
$U^{*-1}$
 beschr\"ankt sind }
\par\addvspace{\smallskipamount}\noindent  
{\rm 2)}
$
 S_+(c^{-1}U)=S_+(c^{-1}U^{*-1})=S_+(c^{-1}U^{-1})=S_+(c^{-1}U^*)=\{0\}  
 \,.
$
\par\addvspace{\medskipamount}\noindent  

 Ist solch ein Operator 
$U$
 gegeben, so setzen wir
${\cal U} := {\hat U} = U\oplus U^{*-1}$.

 Dann gilt:
\par\addvspace{\medskipamount}\noindent      
 {\rm a)}\quad 
$
 S_+(c^{-1}{\cal U})  =  S_+(c^{-1}{\cal U}^{-1})  =  \{0\}\,
$
 \hfill\makebox[35ex][l]{ 
 (nach Bemerkung {\bf B1-2} );} 
\par\addvspace{\smallskipamount}\noindent      
 {\rm b)}\quad
$\cal U$ 
 ist 
$J$-unit\"ar und symplektisch 
 \hfill\makebox[35ex][l]{ 
 ( nach Hilfssatz {\bf H1-1});}
\par\addvspace{\smallskipamount}\noindent      
 {\rm c)}\quad 
 Die Punkte (i), (ii), (iii) sind erf\"ullt 
 \hfill\makebox[35ex][l]{ 
 ( nach Bemerkung {\bf B1-1}).}
\par\addvspace{\medskipamount}\noindent     

 Nun haben wir jenen 
$U$
 zu konstruieren.
 Das wollen wir im Hifssatze {\bf H2-1} tun,
 erst aber m\"ussen wir an gewisse Definitionen und Tats\"ache
 der Theorie von sogenannten gewichteten Verschiebungsoperatoren
 erinnern.

\begin{Definition}{ D2-1}
$H_0$
 sei ein beliebiger 
 separabler (reeller oder komplexer) Hilbertraum;
 das Symbol 
$(\,,\,)$ 
 stehe f\"ur das Skalarprodukt von
$H_0$
 und mit 
$\{b_n\}_n$
 bezeichne man eine orthonormierte Basis, deren Elemente durch  
$n=...,-1,0,1,...$
 indexiert werden.

 Es sei  $\{u_n\}_n$ ,$n\in {\bf Z}$ 
 eine zweiseitige Zahlenfolge; dabei wird $u_n \not= 0$,$n\in {\bf Z}$
 angenommen.
 Nun bezeichne 
$U$ 
 einen Verschiebungsoperator
\footnote{ 
 voller Name: gewichteter Verschiebungsoperator, 
 von $\{b_n\}_n$, nach rechts. } 
, der durch die Formel 
$$
 U\,:\,b_n\,\mapsto\,\frac{u_{n+1}}{u_n}\,b_{n+1}    \eqno(*)
$$
 erzeugt wird.
\end{Definition}

\begin{Observation}{ O2-1}

 Es kann vielleicht nicht fehlen, auf die Definition von $U$ einzugehen.
 Man bildet den Operator $U$ so:
 
 Vorerst setzt man die Zuordnungsvorschrift $(*)$ auf ganze lineare H\"ulle
 von den Basiselementen $\{b_n\}_n\in {\bf Z}$ 
 zu einem linearen Operator fort.
 Diese Fortsetzung
 ist eindeutig und erzeugt einen linearen dicht definierten abschlie{\ss}baren
 Operator, der hier mit $U_{min}$ bezeichnet wird.
 Die Abschlie{\ss}ung von $U_{min}$ ist genau der Operator $U$.

 Also, der hierdurch definierte Operator
$U$
 ist 
 abgeschlossen, 
 mindestens dicht definiert, injektiv, besitzt dichtes Bild und die Wirkung von
 $U^N$,$U^{*-N}$,
                     $U^{*N}U^N$,$U^{-N}U^{*-N}$ 
 (f\"ur alle ganzen $N$) wird durch die Formeln
$$
 U^N:b_n\mapsto \frac{u_{n+N}}{u_n}b_{n+N}\,;\,
 U^{*-N}:b_n\mapsto \frac{u_n^*}{u_{n+N}^*}b_{n+N}\,;
$$
$$
 U^{*N}U^N:b_n\mapsto {|\frac{u_{n+N}}{u_n}|}^2 b_n\,;\,
 U^{-N}U^{*-N}:b_n\mapsto {|\frac{u_n}{u_{n+N}}|}^2 b_n\,;
$$
 erzeugt.

 Insbesondere ist $U^N$ beschr\"ankt genau dann, wenn
 die Zahlenfolge \\
$\{|u_{n+N}/u_n|\}_n$ beschr\"ankt ist.
$\Box$
\end{Observation}

\begin{Observation}{ O2-2}
 Die Schar $\{b_n\}_n$ 
 ist eine orthonormierte Basis und es gilt auch
$U^Nb_n\perp U^Nb_m$ f\"ur $n\neq m$. Deshalb ist
$$
 {\|U^Nf\|}^2=\sum_n|(b_n,f)|^2\|U^Nb_n\|^2=
 \sum_n|(b_n,f)|^2{|\frac{u_{n+N}}{u_n}|}^2\,
$$
 f\"ur jede $f\in D_{U^N}$.

 Insbesondere 
$$\|U^Nf\|\geq|(b_n,f)||u_{n+N}/u_n|$$ 
 f\"ur alle ganzen 
$n$. 
 Daraus folgt: 

 Ist $f\in H_0\setminus \{0\}$ und sind  $M,a$ gewisse Zahlen, so da{\ss}\   
$$
 \|U^Nf\|\leq Ma^N 
 \mbox{ f\"ur } 
 N=0,1,2,...\,,
$$
 dann gibt es eine Zahl $M'$ ,
 so da{\ss}\  
$$
 |u_N|\leq M'a^N 
 \mbox{ f\"ur } 
 N=0,1,2,...
$$
 F\"ur
$U^{*-1},U^{-1},U^*$, es gelten die analogischen Implikationen. 
 Wir schreiben alle diese Implikationen ausf\"uhrlich heraus:
$$
\begin{array}{llcccccccccc}
 
  &
 \|U^Nf\|            &    \leq      &   Ma^N\,        & \Rightarrow &&
 \,    |u_N|         &    \leq      &   M^{'}a^{N}    & 
 \,  (N=0,1,2,...)   &  \\

 &
 \|U^{*-N}f\|        &     \leq     &   Ma^N\,        & \Rightarrow &&
 \,    |u_N|^{-1}    &     \leq     &   M^{'}a^{N}    &
 \,  (N=0,1,2,...)   &  \\

 &
 \|U^{-N}f\|         &     \leq     &    Ma^N\,       & \Rightarrow &&
 \,    |u_{-N}|      &     \leq     &    M^{'}a^{N}   &
 \,  (N=0,1,2,...)   &   \\

 &
 \|U^{*N}f\|         &     \leq     &    Ma^N\,       & \Rightarrow &&
 \,    |u_{-N}|^{-1} &     \leq     &    M^{'}a^{N}   &
 \,    (N=0,1,2,...) &     \\
\end{array}
$$
 (Wir haben hier bei den Bezeichnungen das Format
$$
\begin{array}{cccccc}
 \exists f \in H_0\setminus \{0\}, M>0  , a>0 
 \forall N \geq 0 
 &\cdots 
 & \Rightarrow &
 \exists M'>0  \forall N \geq 0 
 & \cdots 
\end{array}
$$
gemeint.)
 Am n\"achsten Hilfssatz werden wir n\"amlich solche Folgerungen 
 dieser obigen Implikationen anwenden:

 Es stehe $c$ f\"ur eine reelle positive Zahl. Dann gilt:

$$
\begin{array}{lccccccc}

S_+(c^{-1}U) \not= \{0\}         & \Rightarrow & \exists M' > 0 \forall N \geq 0&
 \,    |u_N|         &    \leq   &   M^{'}(c+1)^{N}    & 
\\

S_+(c^{-1}U^{*-1}) \not= \{0\}   & \Rightarrow & \exists M' > 0 \forall N \geq 0&
 \,    |u_N|^{-1}    &     \leq  &   M^{'}(c+1)^{N}    &
\\

S_+(c^{-1}U^{-1}) \not= \{0\}    & \Rightarrow & \exists M' > 0 \forall N \geq 0&
 \,    |u_{-N}|      &     \leq  &    M^{'}(c+1)^{N}   &
\\

S_+(c^{-1}U^*) \not= \{0\}       & \Rightarrow & \exists M' > 0 \forall N \geq 0&
 \,    |u_{-N}|^{-1} &     \leq  &    M^{'}(c+1)^{N}   &
\\

\end{array}
$$
                                                                $\Box $
\end{Observation}

\begin{Hilfssatz}{ H2-1}
 Man setze
$$
 u_n\,:=\,(c+2)^{|n|\sin(\frac{\pi}{2}\log_2(1+|n|))}\qquad(n=...,-1,0,1,...) ,
$$
 hierbei sei $c>0$ beliebig ausgew\"ahlt.

 Dann ist der entsprechende Verschiebungsoperator  $U$ mit seinem Inverse 
 beschr\"ankt, und
$$
 S_+(c^{-1}U)  =  S_+(c^{-1}U^{*-1})  =  S_+(c^{-1}U^{-1})
   =  S_+(c^{-1}U^*)  =  \{0\}
$$
\end{Hilfssatz}

\begin{Beweis}{ }
 Die Ableitung  der zahlenwertigen Funktion
$$
 x\,\mapsto\,|x|\sin(\frac{\pi}{2}\log_2(1+|x|))
$$
 ist gleich
$$
 (\sin(\frac{\pi}{2}\log_2(1+|x|))+\frac{\pi}{2\ln 2}\frac{|x|}{1+|x|}
 \cos(\frac{\pi}{2}\log_2(1+|x|))){\rm sgn}\,x
$$
 und ihr Absolutbetrag \"ubersteigt der Gr\"o{\ss}e $\alpha:=1+\pi/(2\ln 2)$ 
 nicht.
 Nach dem Mittelwertsatz von Lagrange ist 
$$
 (c+2)^{-\alpha}\,\leq\,|u_{n+1}/u_n|\,\leq\,(c+2)^{\alpha}   \,.
$$
 Somit sind die Operatoren $U$ und $U^{-1}$ beschr\"ankt.

 Nun w\"ahlen wir zwei Zahlenfolge so aus:
$$
 n_k\,:=\,2^{1+4k}-1;\quad m_k\,:=2^{3+4k}-1 \qquad (k=1,2,...)  \,.
$$
 Dann gilt: $n_k,\, m_k \in {\bf N}$,\ 
$n_k\,\to\,+\infty$,$m_k\,\to\,+\infty$
 (f\"ur $k\,\to\,+\infty$), und gleichzeitig 
$$
 u_{n_k}  =  u_{-n_k}  =  (c+2)^{n_k};\quad u_{m_k}^{-1}  =  
 u_{-m_k}^{-1}  =  (c+2)^{m_k}   \,.
$$

 Man sieht, da{\ss} man keine Absch\"atzung von der Art 

$$\begin{array}{clc}
 |u_N|\,\leq\,M'(c+1)^N,\,        &
 |u_{-N}|\,\leq\,M'(c+1)^N,\,     &
 |u_{N}|^{-1}\,\leq\,M'(c+1)^N,\,  \\
                                &
 |u_{-N}|^{-1}\,\leq\,M'(c+1)^N  &
\end{array}$$
 (f\"ur $N=0,1,\cdots$)
 verwirklichen kann. 
 Wendet man hier die
   Observation {\bf O2-2} an, so bekommt man, da{\ss} 
$$
S_+(c^{-1}U)=\{0\},
S_+(c^{-1}U^{-1})=\{0\},  
S_+(c^{-1}U^{*-1})=\{0\}, 
S_+(c^{-1}U^*)=\{0\}; 
$$ 
 was zu beweisen war.

\end{Beweis}
 Der Beweis des Hilfssatzes {\bf H2-1} und mithin 
 des Satzes {\bf S2-1} ist beendet.
\end{Beweis}


\newpage
\section
{ Andere Beispiele von $J$-unit\"aren Operatoren }

 In dieser Sektion zeigen wir noch zwei Operatoren, deren 
 Eigenschaften 
 irgendwie absonderlich aussehen.
 Die allgemeine Elemente von der Konstruktionen sind dieselbe, 
 die wir in den vorigen Sektionen eingef\"uhrt haben:

\smallskip 

$H_0$, 
 das steht f\"ur ein beliebiger separabler Hilbertraum;
$\{b_n\}_n$,
 das steht f\"ur eine orthonormierte Basis von 
$H_0$, die Elemente von dieser Basis werden durch  
$n=...,-1,0,1,...$
 indexiert; au{\ss}erdem  
$$
 H:=H_0\oplus H_0 \,,\quad  
 J(x\oplus y):=y\oplus x \,,\quad 
 {\cal J}(x\oplus y):= -y\oplus x \,,\quad (x,y \in H_0) 
$$ 
 und ist 
$ T: H_0 \to H_0$ 
 ein linearer Operator, 
 dann setzen wir
$ \hat T:= T\oplus T^{* -1}$, 
 falls aber 
$T^{* -1}$
 existiert.  

\smallskip 
 
 Wir konstruieren zwei zweiseitige Zahlenfolgen 
$\{v_n\}_n$, 
$\{w_n\}_n$, 
$n=...,-1,0,1,...$, 
 so da{\ss} die entsprechenden Verschiebungsoperatoren, 
$V$ und $W$,  
 und die $J$-unit\"aren und symplektischen Operatoren, 
$\hat V$ und $\hat W$, 
 bemerkliche Eigenshaften haben.

\begin{Definition}{ D3-1} 
 Es sei eine Zahl
$c\geq 1$
 beliebig ausgew\"ahlt.
 Es sei
$v_n:= (2c)^{-|n|}$ 
 f\"ur beliebige ganze 
$n$.
$V :H \to H$ 
 sei der entsprechende,  
 durch die Formeln
$$
 V\, : \, b_n \mapsto \frac{v_{n+1}}{v_n} b_{n+1} 
$$ 
 erzeugte Verschiebungsoperator. 
 Mit anderen Worten, es sei 
$$
 Vb_n :=\frac{1}{2c} b_{n+1} \, \hbox{ f\"ur }   n=0,1,2,... \quad
 Vb_n := 2cb_{n+1} \, \hbox{ f\"ur }     \quad n=...,-2,-1.
$$ 
\end{Definition}

\begin{Bemerkung}{ B3-1}
 Der hierdurch definierte Operator 
$V$ ist beschr\"ankt und 
 besitzt ein beschr\"anktes Inverses;
 nach seiner Definition kann man zeigen, da{\ss} 
\par\addvspace{\smallskipamount}\noindent     
 {\rm (1)} \quad   
$
 \| V^N b_n \|   =   (2c)^{-|n+N|+|n|} \quad
    \mbox{ f\"ur beliebige ganze} \quad n, N;
$ 
\par\addvspace{\smallskipamount}\noindent     
 {\rm (2)} \quad 
$
 r(V) =  r (V^{-1}) = 2c;
$ 
\par\addvspace{\smallskipamount}\noindent     
 {\rm (3)} \quad  
$
  \overline{ S_0 \left(\frac{3c}{2} V \right)} = H, \,
    S_0 \left(\frac{2}{3c} V^{*-1} \right) = \{0\}, \,   
  \overline{ S_0 \left(\frac{3c}{2} V^{-1} \right)} = H, \,
    S_0 \left(\frac{2}{3c} V^{*} \right) = \{0\}.   
$ 
\\
\end{Bemerkung}

\begin{Hilfssatz}{ H3-1}
$L$, $M$ 
 seien solche (in  $\hat H$ lineare abgeschlossene) Teilr\"aume, da{\ss}  
$$
 \hat V L   =   L , 
 \quad | spectrum\, \hat V | L | \, \le \, c \, ,
 \hat V^{-1} M   =   M , 
 \quad | spectrum\, \hat V^{-1} | M | \, \le \, c \, .
$$
 Dann gilt:     
\par\addvspace{\smallskipamount}\noindent     
{\rm (a)}\hfill 
\parbox[t]{.93\textwidth}{ 
$  L   =   L_1 \oplus \{ 0\}, \,
  \quad M   =   M_1 \oplus \{ 0\}, \,
  \mbox{ f\"ur gewisse } \,
       L_1 \, \subset \, H;\, 
       M_1 \, \subset \, H ;
$
}
\par\addvspace{\smallskipamount}\noindent     
{\rm (b)}\hfill 
\parbox[t]{.93\textwidth}{ 
$  VL_1   =   L_1$, $| spectrum\, V | L_1 | \, \le \, c;$
\\ 
$  V^{-1}M_1   =   M_1$, $| spectrum\, V^{-1} | M_1 | \, \le \, c;$
}
\par\addvspace{\smallskipamount}\noindent     
{\rm (c)}\hfill  
\parbox[t]{.93\textwidth}{  
$  L_1 \, \ne \, H,\quad  M_1 \, \ne \, H.$ 
}
\end{Hilfssatz}

\begin{Beweis}{ }

 Beweis von (a) :
 Dies folgt daraus, da{\ss}
\begin{eqnarray*} 
  L \, \subset \, S_0 \left( \frac{2}{3c} \hat V \right)
&=&
  S_0 \left( \frac{2}{3c} V \oplus \frac{2}{3c} V^{*-1} \right)
\\
&=&
  S_0 \left( \frac{2}{3c} V \right) \, \oplus \, 
                    S_0 \left( \frac{2}{3c} V^{*-1} \right)   =  
                    S_0 \left( \frac{2}{3c} V \right) \, \oplus \, \{ 0 \};
\\
 M \, \subset \, S_0 \left( \frac{2}{3c} \hat V^{-1} \right)
&=&
  S_0 \left( \frac{2}{3c} V^{-1} \oplus \frac{2}{3c} V^{*} \right)
\\
&=&
  S_0 \left( \frac{2}{3c} V^{-1} \right) \, \oplus \, 
                    S_0 \left( \frac{2}{3c} V^{*} \right)   =  
                    S_0 \left( \frac{2}{3c} V^{-1} \right) \, \oplus \, \{ 0 \};
\end{eqnarray*}

 Beweis von (b): 
 Nach dem Beweis von (a) haben wir:
$$
 L = L_1 \oplus \{ 0 \}, \, M = M_1 \oplus \{ 0 \}, \, \hat V=V\oplus V^{*-1}.
$$
 Folglich 
$$ 
 spectrum\,\hat V|L = spectrum\,V|L_1,  \quad  
 spectrum\,\hat V^{-1}|M = spectrum\, V^{-1}|M_1, 
$$
 und demnach 
$$ |spectrum\,V|L_1| =| spectrum\,\hat V|L| \leq c,  $$ 
$$ |spectrum\,V^{-1}|M_1| = |spectrum\,\hat V^{-1}|M| \leq c. $$

 Beweis von (c) : 
 Wir haben   
\\
$ | spectrum\, V | L_1 | \, \le \, c$,  
 und
$r(V)=2c$.
 Folglich 
$ L_1  \ne  H$.
 Analog, 
$ | spectrum\, V^{-1} | M_1 | \, \le \, c$,   
 und
$r(V^{-1})=2c$ . 
 Folglich 
$ M_1  \ne  H.$ 

\end{Beweis}

 Nun erinnern wir daran, da{\ss} 
$H\oplus \{0\}$  ein in $\hat H$ 
$J$-neutraler Teilraum ist (sieh [Krein65]). 
Insbesondere ist $H\oplus \{0\}$ ein semidefiniter Teilraum.

 Dann erhaltet man den folgenden

\begin{Satz}{ S3-1}
 Es sei $L$ ein solcher semidefiniter Teilraum von 
$\hat H$, 
 da{\ss}  
$$
\hat V L   =   L, \qquad |spectrum\, \hat V |L| \, \le \, c .
$$
 Dann ist $L$ nicht maximal.

 Es sei $M$ ein solcher semidefiniter Teilraum von 
$\hat H$, 
 da{\ss}  
$$
\hat V^{-1} M   =   M, \qquad |spectrum\, \hat V^{-1} |M| \, \le \, c .
$$
 Dann ist $M$ nicht maximal.
\end{Satz}

\begin{Bemerkung}{ B3-2}

 Der Operator $V$ besitzt die folgende interessante Eigenschaft:

$b_n \in S_0(V) \cap S_0(V^{-1})$ 
 f\"ur jede ganze $n$; somit ist 
$S_0(V) \cap  S_0(V^{-1})$ 
 in $H$ dicht. Noch mehr, $L_k$ bezeichne abgesclossene
 lineare H\"ulle der Menge 
$\{V^s b_k|\, s\ge k\}$. 
 Dann ist  
$$
 VL_k\,\subset L_k \ ,
 r(V|L_k) \,\le\,1/2c \ 
 \mbox{ und } 
 H = \overline{\cup\{L_k\,|k=...-1,0,1...\} }\ .
$$
 aber 
$r(V)=2c$ .

 Der Operator  $V^{-1}$ besitzt die \"ahnliche Eigenschaft.
 Hier ist $L_k$ aber durch abgeschlossene lineare H\"ulle von
$\{V^{-s}b_k|\,s\ge k \}$ 
 zu ersetzen.

 Wir wollen noch eine Eigenschaft von $V$ anmerken.
 
 Es sei 
$f:= {\sum}_{n\not=0}|n|^{-1} b_n$. 
 Dann 
$f \in H_0$
 und 
\begin{eqnarray*}
 {\|V^Nf\|}^2
 &=& \sum_n|(b_n,f)|^2{|\frac{v_{n+N}}{v_n}|}^2\,
\\
 &=&\sum_{n>0}\frac{1}{|n|^{2}}\frac{(2c)^{-2|n+N|}}{(2c)^{-2|n|}}\,
    +\sum_{n<0}\frac{1}{|n|^{2}}\frac{(2c)^{-2|n+N|}}{(2c)^{-2|n|}}\,
\\
 &=&\sum_{n>0}\frac{1}{|n|^{2}}\frac{(2c)^{-2|n+N|}}{(2c)^{-2|n|}}\,
    +\sum_{n>0}\frac{1}{|n|^{2}}\frac{(2c)^{-2|n-N|}}{(2c)^{-2|n|}}\,
\end{eqnarray*}

 Wir wollen die Folge 
$ {\|V^Nf\|}^2 $ 
 {\bf nach unten } absch\"atzen.
 
 Zuerst sei 
$N\geq 0$. 
 Dann haben wir 
\begin{eqnarray*}
 {\|V^Nf\|}^2 
 &\geq&
    +\sum_{n>N}\frac{1}{|n|^{2}}\frac{(2c)^{-2|n-N|}}{(2c)^{-2|n|}} 
\\
 &&
    =\sum_{n>N}\frac{1}{|n|^{2}}\frac{(2c)^{-2n+2N}}{(2c)^{-2n}} 
\\
 &&
    =\sum_{n>N}\frac{1}{|n|^{2}}(2c)^{2N} 
\\
 &&
    \geq\frac{1}{N+1}(2c)^{2N}
    \geq\frac{1}{2^N}(2c)^{2N}
    = 2^N c^{2N} =2^{|N|} c^{2|N|}  \,.
\end{eqnarray*}

 \"Uberdies sehen wir, da{\ss} f\"ur die jetzt zu betrachtenden 
$V$ und $f$ 
 die Gr\"o{\ss}e  
${\|V^Nf\|}^2$
 von  
$N$ 
 so abh\"angt, da{\ss} 
${\|V^Nf\|}^2 = {\|V^{-N}f\|}^2$ 
 f\"ur alle  
$N$
 ist.  

 Somit sehen wir, da{\ss} 
$$ 
 {\|V^Nf\|}^2 
    \geq\frac{1}{|N|+1}(2c)^{2N}
       \geq 2^{|N|} c^{2|N|} \mbox{ f\"ur alle } N \,. 
$$ 
 Ein sehr schnelles Wachstum von  
$ {\|V^Nf\|}^2 $, 
 wenn  
$N \to \pm \infty $
 !!

\end{Bemerkung}

 Jetzt zeigen wir noch einen Beispiel von 
$J$-unit\"arem und symplektischem Operator.

\begin{Definition}{ D3-4}
 Es sei
$w_n:= 2^{-|n|} =2^n$ f\"ur $n \le 0$ und $ w_n:= 1/(n+1)$ f\"ur $n > 0$.
$W :H \to H$ 
 sei der entsprechende,  
 durch die Formeln
$$
 W\, : \, b_n \mapsto \frac{w_{n+1}}{w_n} b_{n+1} 
$$ 
 erzeugte Verschiebungsoperator
 
\end{Definition}

\begin{Bemerkung}{ B3-3}
 Da  $1/2 \le w_{n+1}/w_n \le 2$ 
 f\"ur alle Ganzen $n$ ist, ist der Operator $W$ beschr\"ankt, invertierbar
 und $W^{-1}$ ist beschr\"ankt auch.
\end{Bemerkung}

\begin{Hilfssatz}{ H3-2}
 Der soeben definierte Operator  $W$ hat die Eigenschaften:

$$
 \overline{S_0 (W)}  =   H, \quad 1 \le |spectrum\, W | \le 2, \quad r(W)= 2,
$$
$$
 S (W^{*-1})   =   \{ 0\}, \quad \frac12 \le |spectrum\, W^{*-1} | \le 1, 
\quad r (W^{*-1}) = 1   \,.
$$
\end{Hilfssatz}

\begin{Beweis}{. }
 Der Beweis st\"utzt sich auf die wohlbekannte Spektralradiusformel,
 auf die Bemerkung {\bf B3-3} und auf die Formeln der Observation {\bf O2-1}. 
 Wir haben: 
$$
 \|W^N\|  =   sup \{ \frac{w_{n+N}}{w_n} \, | \, n= \dots -1, 0, 1, \dots \}, 
$$
$$
 \| W^{*-N} \|   =   sup \{ \frac{w_n}{w_{n+N}} \, | \, 
 n= \dots -1, 0, 1, \dots \} .
$$

 Nehmen wir $N>0$ an, und betrachten 
 den Ausdruck $w_n /w_{n+N}$
 im einzelnen.
 Wir sehen: 


\begin{description}
\item[\rm a)]
$ w_n /w_{n+N}   =   1/2^N $  f\"ur $n +N \le 0$;\\
\item[\rm b)]
$ w_n /w_{n+N}   =   (1+n+N)/(1+n)   =   N/(n+1)\, + \, 1 
\, \le \,  N+1$ \hspace*{7mm} f\"ur $0 < n$;\\
\item[\rm c)]
$ w_n /w_{n+N}   =   2^n (1+n+N)\, \le \, N+1$ f\"ur $n\le 0 < N+n$
\end{description}

 Au{\ss}erdem ist $w_0/w_N  =  1+N$ . Folglich ist $ \|W^{*-N}\|\,=N+1$
 (f\"ur $N>0$), daraus $r(W^{*-1})  =  1$ und $ r(W^{-1})=1$ . Ganz analog  
 ergibt sich, da{\ss} $\|W^N\|  =  2^N$ (f\"ur $N>0$), daraus
$r(W)   =  2$ und $ r(W^*)  =  2$.

 Zuletzt, ist $n+N  >  0$,
 so ist $W^N b_n   =  {w_n}^{-1}(1+N+n)^{-1}b_{n+N}$.
 Daraus folgt, da{\ss}  $b_n\in S_0(W)$ f\"ur alle ganzen Zahlen $n$.
 Folglich ist
$\overline{S_0(W)}  =  H$  und \mbox{$ S(W^{*-1})=\{0\}\,. $}

\end{Beweis}

\begin{Satz}{ S3-2}

 Es gibt einen $J$-unit\"aren Operator $\hat W$ 
 und einen maximalen semidefiniten Teilraum   
$L$, so da{\ss} gilt:  
\par\addvspace{\bigskipamount}\noindent     
{\rm (a)}\hfill 
\parbox[t]{.93\textwidth}{ 
 $\hat WL^\perp =  L^\perp 
 \,,\quad 1\le |spectrum\hat W|L^\perp \,| \, \le \, 2,\quad
   r(\hat W|L^\perp )\,=2 $ 
 \medskip \\ 
 trotzdem ist 
$L^\perp = \overline{S_0(\hat W)} \,=\, \overline{S(\hat W)}$.\\
 }
\par\addvspace{\smallskipamount}\noindent     
{\rm (b)}\hfill 
\parbox[t]{.93\textwidth}{ 
$ \hat WL =  L ,\quad |spectrum\hat W|L | \, \le \, 1 $,
 \medskip \\ 
 jedoch ist $ L \cap \overline{ S(\hat W )} = \{0\}$. 
 }

\end{Satz}

\begin{Beweis}{ }
 Man setze 
$$
 L:=\{0\}\oplus H \,,
 \, M:=H\oplus \{0\} \equiv L^\perp \,,
$$ 
 und wende 
 den nun bewiesene Hilfssatz {\bf H3-2} 
 auf die Formeln f\"ur $S_0(\hat W)$
 und $S(\hat W)$ (sieh Einleitung) an:
$$
 S_0(\hat W)\,=S_0(W)\oplus S_0(W^{*-1})\,=S_0(W)\,\oplus \{0\}\,\subset M 
$$
$$
 S(\hat W)\,=S(W)\oplus S(W^{*-1})\,=S(W)\,\oplus \{0\}\,\subset M
$$
 Die Mengen $S_0(W)$ und $S(W)$ beide sind aber in $H$ dicht.
 Folglich fallen 
 ihre Abschlie{\ss}ungen mit $M$ zusammen.
 Um Beweis zu vollenden, ist es hinreichend  anzumerken:
 der Operator 
$\hat W|M $ 
 ist zu 
$W$ 
 unit\"ar \"aquivalent,
 der Operator 
$\hat W|L $ 
 ist zu 
$W^{*-1}$
 unit\"ar \"aquivalent,
 und danach nochmals den Hilfssatz {\bf H3-2} anzuwenden.
\end{Beweis}

\newpage


\section
{ Ein \"Ubergang zu einem kontinuierlichen Modell }

 Dieser \"Ubergang wird mit solchen wohl \"ublichen Schritten erf\"ullt:

$$
 V^N: b_n \mapsto \frac{v_{n+N}}{v_n} b_{n+N}
$$
$$
 V^N: \sum_n f(n) b_n \mapsto \sum_n \frac{v_{n+N}}{v_n}f(n)b_{n+N}
                            =\sum_n \frac{v_n}{v_{n-N}}f(n-N)b_n\,;
$$
$$
 V^N: f(n) \mapsto \frac{v_n}{v_{n-N}}f(n-N)
$$
$$
 V(t): f(x) \mapsto \frac{v(x)}{v(x-t)}f(x-t)
$$
 Dabei ergibt sich, da\ss
$$
\begin{array}{rcr}
 V(t)V(\tau )^{-1} \, : 
 &&
 \\
 f(x) 
 & 
 \stackrel{V(\tau )^{-1}}{\longmapsto}\dfrac{v(x)}{v(x+\tau)}f(x+\tau)
      \stackrel{V(t)}{\longmapsto}
 & \dfrac{v(x)}{v(x-t)}\dfrac{v(x-t)}{v(x-t+\tau)}f(x-t+\tau) %
 \\&&{} =\left( V(t-\tau)f \right)(x) 
\end{array}
$$
 Mit anderen Worten ist die durch $V$ erzeugte Dynamik Zeit-autonom
 und formaler Generator sieht so aus:
\begin{eqnarray*}
\left(Hf\right)(x)&=&\left(V(t)f\right)_{t=0}'(x)  \\
&=&\frac{\partial}{\partial t}\left[\frac{v(x)}{v(x-t)}f(x-t)\right] _{t=0} \\
&=&-\frac{\partial f(x)}{\partial x} +\frac{v'(x)}{v(x)}f(x)
  = -v(x)\frac{\partial }{\partial x}\Bigl(\frac{1}{v(x)}f(x)\Bigr)  \,.
\end{eqnarray*}

 Wenden wir insbesondere den hier geschriebenen \"Ubergang
 auf diejenigen diskreten Systeme an,   
 die wir in vorhergehenden 
 Sektionen betrachtet haben.
 Dann lassen sich
 die entsprechenden kontinuierlichen Systeme 
 so beschreiben:
\begin{description}
\item[\rm a)]
$v(x)=e^{|x|sin\big(ln(1+|x|)\big)}$ \\
$(Hf)(x)$  \\
${ }=-\dfrac{\partial f(x)}{\partial x}
 +\Big(sin(ln(1+|x|))
  +\dfrac{|x|}{1+|x|}cos\bigl(ln(1+|x|)\bigr)\Big)sgn(x)f(x) $
\\
\item[\rm b)]
$ v(x)=e^{-|x|}$\\
$(Hf)(x)=-\dfrac{\partial f(x)}{\partial x}-(sgn\,x)f(x)$\\
\item[\rm c)]
$$
 v(x)=\left\{\begin{array}{lcl} 
 e^x&,&x<0\\
\dfrac{1}{x+1}&,&x>0 
\end{array} \right\}
$$
$$
 (Hf)(x)=-\frac{\partial f(x)}{\partial x}
 +\left\{
 {1 ,x<0 \atop -\dfrac{1}{x+1} ,x>0 }
 \right\}
 f(x) 
$$
\end{description}
 
\begin{Bemerkung}{.}
 Alle  vorstehenden Beispiele 
 von $J$-unit\"aren Operatoren (symplektischen Automorphismen)
 sind durch 
 Hilbertsche direkte Summen von gewichteten zweiseitigen 
 Verschiebungsoperatoren realisiert.

 Somit werden diese Operatoren unit\"ar sein,
 wenn man den zugrundeliegenden Hilbertraum geeignet (aber un\"aqivalent !)
 renormiert.

 Also, die Operatoren sind `gut' und hingegen sind die R\"aume
 `schlecht'.

 Ich hoffe und glaube, da{\ss} dies typisch ist.
 Jetzt aber kann ich einige S\"atze nur von solcher Art beweisen:

1)
 ${\cal H}$ 
 sei ein linearer Raum, der ist mit einem inneren Produkt
\footnote{ indefiniten oder definiten, das hat hier keine Bedeutung }
 \mbox{$<\cdot,\cdot>$},
 oder mit einer symplektischen bilinearen Form 
 \mbox{$s(\cdot,\cdot)$} 
 versehen. 
 Ein linearer Operator 
${\cal T}: {\cal H}\to {\cal H}$ 
 sei
 \mbox{$<\cdot,\cdot>$}-unit\"ar (oder symplektisch). 
 Nehmen wir an: es gibt ein Element 
 $u \in {\cal H}$, 
 so da{\ss} 
 die Elemente 
$\{{\cal T}^N u\}_N$
$(N = ...-2,-1,0,1,2...)$ 
 linear unabh\"angig sind. 

 Dann gibt es einen zweiseitigen gewichteten Verschiebungsoperator 
$U: H_0 \to H_0$
 und einen 
$f\in H_0\oplus H_0$, 
 so da{\ss} jeweils gilt: 
$$
 <u, {\cal T}^N u> = (f,J\hat U^{N}f) 
 \qquad (N=...-2,-1,0,1,2,...) 
$$ 
 oder 
$$
 s(u, {\cal T}^N u) = (f,{\cal J}\hat U^{N}f) 
 \qquad (N=...-2,-1,0,1,2,...) 
$$

2)
${\cal T}$ sei beliebiger invertierbarer Operator, angenommen aber, da{\ss} 
$$
 S_0({\cal T}) = \{0\}\,,\,
 S_0({\cal T}^{*-1}) = \{0\}\,.
$$
 Dann ist ${\cal T}$ ein s-Limes von eine Folge von Operatoren,
 die zu unit\"aren Operatoren \"ahnlich sind [Ch00].
\end{Bemerkung}

\newpage 





\newpage
\bibliographystyle{unsrt}

\end{document}